\documentclass[letterpaper, 10 pt, conference]{ieeeconf}  

\IEEEoverridecommandlockouts                              

\IEEEoverridecommandlockouts

\usepackage{cite}
\usepackage{amsmath,amssymb,amsfonts}
\usepackage{graphicx}
\usepackage{textcomp}
\usepackage{xcolor}
\usepackage{algpseudocode}

\usepackage{amsthm}
\newtheorem{lemma}{Lemma}
\newtheorem{proposition}{Proposition}  

\usepackage[linesnumbered,ruled,vlined]{algorithm2e}
\def\BibTeX{{\rm B\kern-.05em{\sc i\kern-.025em b}\kern-.08em
    T\kern-.1667em\lower.7ex\hbox{E}\kern-.125emX}}
\begin{document}

\title{Multi-armed Bandit for Stochastic Shortest Path in Mixed Autonomy\\
}


\author{Yu Bai, Yiming Li, and Xi Xiong
\thanks{This work was supported in part by NSFC Project 72371172,  NSFC Project 52272320 and Fundamental Research Funds for the Central Universities.}
\thanks{Y. Bai, Y. Li and X. Xiong are with the Key Laboratory of Road and Traffic Engineering, Ministry of Education, Tongji University, Shanghai, China (Emails: baiyu@tongji.edu.cn, 2410792@tongji.edu.cn, xi\_xiong@tongji.edu.cn,)}
}

\maketitle

\begin{abstract}
In mixed-autonomy traffic networks, autonomous vehicles (AVs) are required to make sequential routing decisions under uncertainty caused by dynamic and heterogeneous interactions with human-driven vehicles (HDVs). Early-stage greedy decisions made by AVs during interactions with the environment often result in insufficient exploration, leading to failures in discovering globally optimal strategies. The exploration–exploitation balancing mechanism inherent in multi-armed bandit (MAB) methods is well-suited for addressing such problems. Based on the Real-Time Dynamic Programming (RTDP) framework, we introduce the Upper Confidence Bound (UCB) exploration strategy from the MAB paradigm and propose a novel algorithm. We establish the path-level regret upper bound under the RTDP framework, which guarantees the worst-case convergence of the proposed algorithm. Extensive numerical experiments conducted on a real-world local road network in Shanghai demonstrate that the proposed algorithm effectively overcomes the failure of standard RTDP to converge to the optimal policy under highly stochastic environments. Moreover, compared to the standard Value Iteration (VI) framework, the RTDP-based framework demonstrates superior computational efficiency. Our results highlight the effectiveness of the proposed algorithm in routing within large-scale stochastic mixed-autonomy environments.
\end{abstract}

\section{INTRODUCTION}

With the development of autonomous driving technologies, autonomous vehicles (AVs) and human-driven vehicles (HDVs) increasingly share road space in mixed-autonomy environments~\cite{li2022cooperative}. In such settings, AVs are required to make real-time routing decisions under the behavioral variability and uncertainty of HDVs, posing significant planning challenges. Although substantial progress has been made in AV-HDV interaction modeling and individual vehicle control~\cite{wang2022real}, global routing for AVs in mixed-autonomy environments remains underexplored and calls for further investigation. AVs equipped with dynamic routing capabilities can proactively avoid potential congestion in mixed-autonomy environments~\cite{Sun201532}, thereby alleviating traffic pressure and enhancing the overall efficiency and stability of the transportation system. 
\begin{figure}[h!]
    \centering
\includegraphics[width=0.9\columnwidth]{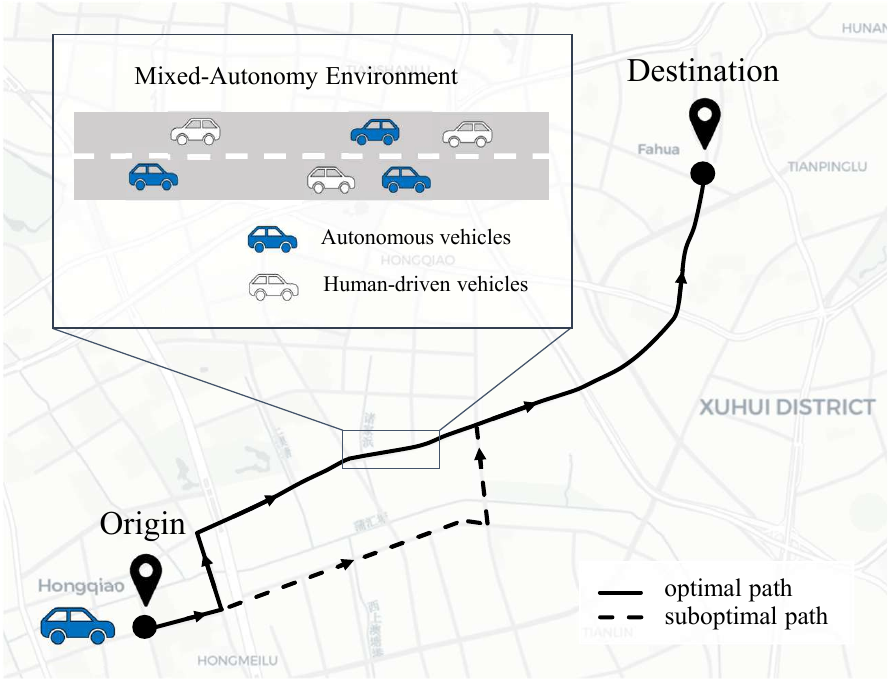}  
    \caption{AVs Routing Scenario in Mixed-Autonomy Environments}
    \label{fig:scenario }
\end{figure}

In this study, we focus on the path routing problem for AVs in mixed-autonomy environments with dynamically changing road conditions. As illustrated in Figure~\ref{fig:scenario }, we construct a local urban network shared by AVs and HDVs, where edge weights and state transitions are subject to stochastic perturbations induced by human driving behavior. To address this, we propose a path routing algorithm that integrates Real-Time Dynamic Programming (RTDP) with Upper Confidence Bound (UCB)-based exploration. The algorithm selects paths based on the principle of optimism and accelerates convergence through asynchronous dynamic updates, enabling AVs to efficiently navigate highly stochastic traffic networks and progressively approach near-optimal routing policies.

Based on the Markov Decision Process (MDP) framework, Bertsekas systematically defined the Stochastic Shortest Path (SSP) problem~\cite{bertsekas2012dynamic}, providing a theoretical foundation for modeling the randomness in mixed-autonomy environments ~\cite{pattanamekar2003dynamic}. From a model-based perspective, Guillot and Stauffer leverage value iteration and policy iteration methods to systematically study the SSP problem~\cite{guillot2020stochastic}. From a model-free perspective, reinforcement learning methods have been widely applied. Hoel et al. proposed a lane-change decision-making method based on deep reinforcement learning, which enables autonomous learning of optimal policies through interaction with the environment~\cite{Hoel20182148}. However, it often suffers from slow convergence when dealing with large-scale state spaces. Barto et al. proposed the Real-Time Dynamic Programming (RTDP) algorithm, which, based on the asynchronous DP theory, significantly improves convergence efficiency in large-scale state spaces~\cite{barto1995learning}. However, the greedy strategy of RTDP may miss globally optimal policies in highly stochastic environments. Methods inspired by the Multi-Armed Bandit (MAB) framework offer promising directions for improving exploration strategies~\cite{lattimore2020bandit,8895728}. Tarbouriech et al. introduced an exploration bonus mechanism within the Value Iteration framework~\cite{tarbouriech2021stochastic}, while Zhang et al. extended the UCB-based exploration strategy to graph structures~\cite{zhang2023multi}, aiming to alleviate the exploration deficiency present in traditional methods. 

Therefore, this paper integrates the UCB exploration mechanism from the MAB framework into RTDP, proposing a novel algorithm for solving SSP problems. In each interaction, the algorithm estimates the costs of feasible paths based on the principle of optimism and guides exploration by selecting the path with the currently lowest estimated cost. Moreover, it adopts an asynchronous update scheme that incrementally updates only the state values along the traversed path, thereby avoiding the high computational cost of global updates. Furthermore, we derive a path-level regret bound under the RTDP framework, offering theoretical insights into the convergence behavior and performance guarantees of the proposed algorithm. Importantly, our algorithm targets efficient exploration and fast convergence under high stochasticity. Based on this, the contributions of our work are summarized as follows. 
\begin{itemize}
\item We achieve innovation in the algorithmic structure by integrating RTDP with the MAB framework, and propose a novel algorithm that combines local updates with optimistic exploration strategies for decision-making. This integration effectively addresses the challenge of balancing exploration and exploitation for AVs routing through SSP problems in highly stochastic mixed-autonomy environments. 
\item We establish a novel path-level regret bound of order $\mathcal{O}(\sqrt{T \log T})$ under the RTDP framework, providing theoretical guarantees for convergence and algorithmic reliability. 
\item We evaluate our algorithm on a real-world local road network in Shanghai. Results show that standard RTDP struggles to converge, while our algorithm achieves stable convergence and near-optimal path strategies. Compared to value iteration with exploration (e.g., Value\_Iteration\_UCB), our algorithm also shows superior computational efficiency. 
\end{itemize}

The remainder of this paper is organized as follows. Section II presents the modeling and formulation of mixed-autonomy environments as an SSP problem. Section III details the proposed algorithm, including the algorithmic design, the derivation of the path-level regret bound, and the theoretical analysis of its computational complexity. Section IV evaluates the performance of the proposed algorithm against several benchmark algorithms, based on experiments conducted on a real-world local road network in Shanghai. Finally, Section V concludes the paper by summarizing the main findings and discussing potential future research directions for AV routing problems.

\section{Modeling and Formulation} \label{02_modeling}
We model the routing problem of autonomous vehicles (AVs) in mixed-autonomy environments as a Stochastic Shortest Path (SSP) problem~\cite{zimin2013online}. 
Given a directed graph \( G = (\mathcal{V}, \mathcal{E}) \), where each vertex \( v \in \mathcal{V} \) represents an intersection and each edge \( e \in \mathcal{E} \) corresponds to a road segment, the SSP problem is formulated as a Markov Decision Process (MDP) \( \mathcal{M} = (\mathcal{S}, \mathcal{A}, \mathcal{P}, \mathcal{C}) \)~\cite{bertsekas2012dynamic}.
Here, \( \mathcal{S} \) denotes the finite set of states (i.e., nodes in the road network), \( \mathcal{A}(s) \) is the set of available actions at state \( s \in \mathcal{S} \), each corresponding to choosing an outgoing edge from \( s \).
In our setting, the state transition function \( \mathcal{P}: \mathcal{S} \times \mathcal{A} \rightarrow \mathcal{S} \) is deterministic, mapping each state-action pair \( (s,a) \) to a unique next state \( s' = \mathcal{P}(s,a) \). The cost function \( \mathcal{C} \) maps each state-action pair \( (s,a) \) to a random variable \( C(s,a) \), representing the stochastic cost. A realization of this cost is denoted by \( c(s,a) \sim C(s,a) \), and its expected value is given by \( \mathbb{E}[C(s,a)] \).

The objective is to find a stationary policy \( \pi: \mathcal{S} \rightarrow \mathcal{A} \) that minimizes the expected cumulative cost from any initial state \( s_o \in \mathcal{S} \) to a designated terminal state \( s_d \in \mathcal{S} \). The optimal value function \( V^*(s) \) is defined as the minimum expected total cost required to reach \( s_d \) from state \( s \) under the optimal policy \( \pi^* \), i.e.,
\[
V^*(s) = \mathbb{E}\left[ \sum_{t=0}^{T-1} c(s_t, a_t) \mid s_o = s, \pi^* \right],
\]
where \( T \) is the (random) time at which the terminal state \( s_d \) is reached.

We consider an AV that, in each episode, starts from $s_{o}$, executes actions while paying a certain cost $c(s,a)$, and transitions to new states until reaching the destination $s_{d}$. The path taken during this episode is denoted as $P$. Since the AV initially knows only the graph structure but not the distribution of the stochastic road segment characteristics (i.e., edge weights), it can only estimate them through online learning. During each episode, the AV updates the value function and learns decision-making through real-time interactions with the environment. The expected shortest path is denoted as $P^*$, and the optimal state value function as $V^*(s)$.  After $T$ episodes, we define regret to evaluate the agent's performance as follows: 
\begin{equation}
R_T :=\sum_{t=1}^{T}( \sum_{e \in P_t} \mu_e - \sum_{e \in P^*} \mu_e).
\label{eq:regret}
\end{equation}

To address this online decision-making problem under uncertainty, the agent not only estimates the edge cost distributions through interaction but also balances exploration and exploitation to minimize cumulative regret. This motivates the development of an efficient learning algorithm that integrates exploration strategies into the planning process, which we achieve by combining the UCB principle with the RTDP framework.

\section{Multi-Armed Bandit for adaptive routing in mixed-autonomy} 
\label{03_greedy_algorithm}

In this section, we propose an algorithm that incorporates the UCB strategy into the RTDP framework to solve the routing problem of AVs.  In this algorithm, the Bellman equation is used to update the value function $V(s)$ for each state, ensuring the convergence of the iterative process. Additionally, the UCB decision-making strategy at each state effectively balances exploration and exploitation, preventing a greedy policy from prematurely converging to a suboptimal path. Subsequently, we derive a path-level regret bound, showing sublinear cumulative regret and efficient convergence within a finite number of episodes. 

\subsection*{A. Overall Framework}
This section introduces the overall framework of the proposed algorithm, which balances exploration and exploitation within the RTDP framework by integrating a UCB strategy. At each state, the algorithm computes the expected action value for each state-action pair based on the Bellman equation, as follows: 
\begin{equation}
Q(s, a) = \hat{c}(s, a) + \hat V(s'),
\end{equation}
where $\hat V(s')$ denotes the estimated value of the next state. 

Based on the principle of optimism and Lemma~\ref{lemma:Hoeffding}, the optimistic cost estimate is defined as
\begin{equation}
U(s, a) := \hat{c}(s, a) +\hat V(s')- rad(e),
\end{equation}
where the confidence radius is given by
\begin{equation}
rad(e):=  \sqrt{ \frac{2 \log N(s)}{n(e)} },
\end{equation}
where $N(s)$ and $n(e)$ denote the number of visits to state $s$ and edge $e$, respectively.  

At each state, the agent selects the action with the minimum $U$, executes the action, and samples the corresponding edge cost. Based on the sampled cost, the agent updates the value function estimate of the source state according to Equation ~\eqref{eq:update}. 
\begin{subequations}\label{eq:update}
\begin{align}
a &= \arg \min_{(s,a) \in \mathcal{E}(s)} U(s,a), \label{eq:update-a} \\
V(s) &\leftarrow \min\left(V(s), \hat{c}(s,a) + \hat{V}(s')\right), \label{eq:update-b}
\end{align}
\end{subequations}
where $\mathcal{E}(s) := \{ e \in \mathcal{E} \mid \text{source}(e) = s \}.$

This reflects the asynchronous update mechanism of RTDP, which updates the value function only for the current state visited during each episode, instead of performing global updates across all states. This process is repeated along the sampled path until the terminal state $s_{d}$.
\subsection*{B. Main Theorems}
This section presents a detailed derivation of the path-level regret bound under the RTDP framework. In contrast to existing regret bound analyses based on the standard Value Iteration framework~\cite{tarbouriech2021stochastic, auer2008near}, the asynchronous update nature of RTDP and the confidence-based exploration mechanism of UCB influence the approach and structure of the regret bound analysis, where the regret is defined as in Equation~\eqref{eq:regret}.

To derive an upper bound, we relax Equation~\eqref{eq:regret} by decomposing it into two components: the price of optimism term and the Bellman error term. Such a regret decomposition has been widely used in the literature on exploration algorithms and reinforcement learning ~\cite{auer2008near, lattimore2020bandit, zhang2023multi}. 
\begin{align}
R_T &:= \sum_{t=1}^{T} \left( \sum_{e \in P_t} \mu_e - \sum_{e \in P^*} \mu_e \right) \notag \\
&\leq \sum_{t=1}^T \left[
\underbrace{\mu(P_t) - \hat{\mu}(P_t)}_{\text{Price of Optimism}} +
\underbrace{\sum_{s \in P_t} \left( \hat{V}_t(s) - V^*(s) \right)}_{\text{Bellman Error}}
\right].
\label{eq:decomposition}
\end{align}

The following proposition establishes the regret bound of the proposed algorithm.
\begin{proposition}
Let \(T > 0\) be the total number of episodes. Define \(|\mathcal{E}|\) as the number of edges in the graph,  $\mathcal{S}_{\text{rel}}$ as the set of states that are visited at least once across all episodes, and  $L_{\max}$ as the maximum number of steps in all paths. 

The cumulative regret  $R_T$ of the proposed algorithm satisfies the following upper bound: 
\begin{align}
R(T) \leq\ 
&\underbrace{\mathcal{O}\left( \sqrt{|\mathcal{E}| \cdot T \cdot \log T} \right)}_{\text{Price of Optimism}} \notag \\
+ \, 
&\underbrace{\mathcal{O}\left( \sqrt{|\mathcal{S}_{\text{rel}}| \cdot T \cdot \log T \cdot L_{\max}} \right)}_{\text{Bellman Error}}.
\label{bound}
\end{align}
\label{prop:regret bound}
\end{proposition}
\subsubsection{Price of Optimism}
The first term in the decomposition, the \textit{Price of Optimism}, captures the cumulative error arising from estimation uncertainty. It measures the total difference between the true expected costs and the optimistic estimates used across all episodes. 
\begin{lemma}
(Hoeffding's Inequality) Let \( X_1, X_2, \dots, X_n \) be independent random variables bounded in \([0, 1]\). Let \(\bar{X}_n = \frac{1}{n} \sum_{i=1}^n X_i\) denote their empirical mean. Then, for any \( \delta > 0 \), with probability at least \(1 - \delta\),
\[
\left| \bar{X}_n - \mathbb{E}[X] \right| \leq \sqrt{ \frac{1}{2n} \log \frac{2}{\delta} }.
\]
\label{lemma:Hoeffding}
\end{lemma}
Based on Lemma~\ref{lemma:Hoeffding} ~\cite{min2022learning}and the principle of UCB, the confidence radius $\text{rad}_t(e)$ is used to upper bound the cumulative cost estimation error, referred to as the price of optimism. \begin{align}
\text{Price of Optimism} :=\ 
&\sum_{t=1}^{T} \left( \mu(P_t) - \hat{\mu}(P_t) \right) \notag \\
=\ 
&\sum_{t=1}^{T} \sum_{e \in P_t} \text{rad}_t(e) \notag \\
=\ 
&\sum_{t=1}^{T} \sum_{e \in P_t} \sqrt{ \frac{2 \log N_t(s)}{n_t(e)} } \label{eq:po-step3} \\
=\ 
&\sum_{e \in \mathcal{E}} \sum_{k=1}^{n_e(T)} \sqrt{ \frac{2 \log N_k(s)}{k} }.\label{eq:po-step4}
\end{align}

Reordering the double summation in Equation ~\eqref{eq:po-step3}—taken initially over rounds $t$ and edges \( e \in P_t \) —into the form of  Equation~\eqref{eq:po-step4}, where the sum is grouped by edge $e$ and indexed by its  \( k \)-th visitation, makes the cumulative contribution of each edge explicit. 
\begin{align}
\text{Price of Optimism} 
&\leq \sum_{e \in \mathcal{E}} \sum_{k=1}^{n_e(T)} \sqrt{ \frac{2 \log T}{k} } \notag \\
&= \sum_{e \in \mathcal{E}} \sqrt{2 \log T} \cdot \sum_{k=1}^{n_e(T)} \frac{1}{\sqrt{k}}.
\label{eq:exploration-upperbound}
\end{align}

\begin{lemma}
For any positive integer \( n \), the following inequality holds:
\[
\sum_{k=1}^{n} \frac{1}{\sqrt{k}} \leq 2 \sqrt{n}.
\]
\label{lemma:harmonic-root}
\end{lemma}
By applying Lemma~\ref{lemma:harmonic-root}, we can further bound it as follows. 
\begin{align}
\sum_{k=1}^{n_e(T)} \sqrt{ \frac{2 \log T}{k} }
&\leq \sqrt{2 \log T} \cdot \sum_{k=1}^{n_e(T)} \frac{1}{\sqrt{k}} \notag \\
&\leq 2 \sqrt{2 \log T} \cdot \sqrt{n_e(T)}.
\label{eq:harmonic-bound}
\end{align}
\begin{lemma}
(Cauchy–Schwarz Inequality) For any non-negative sequence \(\{n_e(T)\}_{e \in \mathcal{E}}\), it holds that:
\[
\sum_{e \in \mathcal{E}} \sqrt{n_e(T)} \leq \sqrt{|\mathcal{E}| \cdot \sum_{e \in \mathcal{E}} n_e(T)}.
\]
\label{lemma:cauchy-schwarz}
\end{lemma}
Apply Lemma~\ref{lemma:cauchy-schwarz} to upper bound the square root terms in the summation above:
\begin{equation}
\sum_{e \in \mathcal{E}} \sqrt{n_e(T)} 
\leq \sqrt{|\mathcal{E}| \cdot \sum_{e} n_e(T)}.
\label{eq:cauchy-schwarz}
\end{equation}

Since \( n_e(T) < T \), we have:
\begin{align}
\text{Price of Optimism} 
&\leq 2 \sqrt{2 \log T} \cdot \sum_{e \in \mathcal{E}} \sqrt{n_e(T)} \notag \\
&\leq 2 \sqrt{2 \log T} \cdot \sqrt{|\mathcal{E}| \cdot \sum_{e} n_e(T)}.
\label{eq:exploration-final-bound}
\end{align}
As a result, the \textit{Price of Optimism} is bounded above by:
\begin{equation}
\text{Price of Optimism} = \mathcal{O}\left( \sqrt{|\mathcal{E}| \cdot T \cdot \log T} \right).
\end{equation}
\subsubsection{Bellman Error}  
The second term in the decomposition, the \textit{Bellman Error}, represents the deviation between the estimated and true state values over all sampled paths. Due to the asynchronous update mechanism of RTDP, value function updates are performed only on the states visited along each sampled path. Therefore, the analysis of the Bellman error should be based on the set of states visited in each episode.
\begin{align}
\text{Bellman Error} 
&:= \sum_{t=1}^{T} \sum_{s \in P_t} \left( \hat{V}_t(s) - V^*(s) \right)\notag \\
&= \sum_{s \in \mathcal{S}_{\text{rel}}} \sum_{t: s \in P_t} \delta_t(s) \notag \\
&= \sum_{s \in \mathcal{S}_{\text{rel}}} \sum_{n=1}^{N_T(s)} \delta_n(s).
\end{align}

Similarly, by applying Lemma~\ref{lemma:harmonic-root} and considering the constraint on the path length, we can further upper bound the cumulative Bellman error as follows: 
\begin{align}
\sum_{s \in \mathcal{S}_{\text{rel}}} \sum_{n=1}^{N_T(s)} \delta_n(s) 
&\;\le\; 
2 \sqrt{2\log T}  \cdot \sum_{s} \sqrt{N_T(s)}.
\end{align}

By applying Lemma ~\ref{lemma:cauchy-schwarz} and noting that $|P_t|<L_{max}$, we have: 
\begin{align}
\sum_{s \in \mathcal{S}_{\text{rel}}} \sqrt{N_T(s)} 
&\le 
\sqrt{|\mathcal{S}_{\text{rel}}| \cdot \sum_{s} N_T(s)} \notag \\
&= 
\sqrt{|\mathcal{S}_{\text{rel}}| \cdot \sum_{t=1}^T |P_t|} \notag \\
&\le 
\sqrt{|\mathcal{S}_{\text{rel}}| \cdot T \cdot L_{\max}}.
\end{align}

As a result, the \textit{Bellman Error} is bounded above by:
\begin{align}
\text{Bellman Error} 
\;\le\;
O\left( \sqrt{|\mathcal{S}_{\text{rel}}| \cdot T \cdot \log T\cdot L_{\max}} \right).
\end{align}

This term captures the cumulative error from using estimated instead of optimal state values, highlighting the effect of local updates under partial state coverage. 

Therefore, combining the bounds for the Price of Optimism and the Bellman Error, the overall regret  $R_T$ is bounded above as given in Equation~\eqref{bound}, ensuring that the cumulative regret grows sublinearly with $T$. This result ensures diminishing average regret and algorithmic convergence.  

\subsection*{C. The Algorithm}
Algorithm 1 shows the pseudocode of the proposed algorithm.
\begin{algorithm}[h]
    \SetAlgoLined
    \caption{RTDP-UCB }
    \KwIn{%
        Graph $G = (V,E)$,
        Convergence threshold $\theta$ 
    }
    \tcp{Initialization}
    \ForEach{state \( s \in S \)}{
    $V(s) \leftarrow 0$ ;
    $N(s)\leftarrow 0$
    
         \ForEach{edge \( e = (s,a) \in E \)}{
        $n(e) \leftarrow 0$ ;
        $sum_e \leftarrow 0$ ;
        $\hat{c}(s,a) \leftarrow 0$ 
         }       
    }
    \For{each episode}{
        $s \gets s_{o}$ 
        
        \While{$s \neq s_{d}$}{
            \ForEach{$a \in A(s)$}{
                $Q(s, a) \gets \hat{c}(s,a) +  V(s')$ \;
                $U(s,a) \gets Q(s, a) - \sqrt{ \frac{2 \log N(s)}{n(e)} }$\;
            }
            $a^* \gets \arg\min_{a \in A(s)} U(s,a)$ \;
            Sample ${c}(s, a^*)$ and observe $s'$ \;
            \tcp{Update edge cost estimate}
            $n(e) \gets n(e) + 1$ \;
            $sum_e \gets sum_e + {c}(s, a^*)$ \;
            $\hat{c}(s,a)\gets sum_e /n(e)$ \;
            \tcp{Update state value}
            $V(s) \gets \min_{a \in A(s)} Q(s, a)$ \;
            $s \gets s'$ \;
            
        }
    }
\end{algorithm}

Compared to existing algorithms that integrate dynamic programming (DP) frameworks with exploration strategy~\cite{tarbouriech2021stochastic}, the algorithm proposed in this paper exhibits significant advantages in terms of computational efficiency. The standard Value Iteration algorithm synchronously updates all state values in each iteration, with a per-iteration computational complexity of $\mathcal{O}\left( |\mathcal{S}| \times d_{\max} \right)$.

According to the theory of dynamic programming, to achieve convergence of all state values within a specified precision threshold \( \theta \), approximately $\mathcal{O}\left( log(1/\theta) \right)$ iterations are required. Therefore, the overall time complexity of this algorithm is $\mathcal{O}\left( |\mathcal{S}| \times d_{\max}\times log(1/\theta)\right)$, where $d_{\max}$ denotes the maximum out-degree of the graph.

However, the proposed algorithm is based on the asynchronous DP theory. In each episode, it conducts local updates of the value function along the sampled path. So the time complexity of a single episode is $
\mathcal{O}\left( L_{\max} \times d_{\max}\right)$, which is significantly lower than that of approaches based on standard Value Iteration. 

In summary, the proposed algorithm offers both lower computational complexity and theoretical convergence guarantees, making it particularly suitable for large-scale and sparsely connected state spaces.

\section{Numerical Results} 
\label{04_experiments}
We conducted experiments on a real-world subnetwork located in Xuhui District, Shanghai, which consists of 22 intersections. The weight of each edge was modeled as a Gaussian random variable with the mean corresponding to the actual travel time and a fixed variance of 2, i.e.,$\mu_i \sim \mathcal{N}(t_i, \sigma_i^2)$, to simulate the stochasticity inherent in mixed-autonomy environments.  The proposed algorithm was compared against the following benchmark algorithms:
\begin{itemize}
\item RTDP\_Standard: The standard RTDP algorithm based on the theory of asynchronous dynamic programming\cite{barto1995learning}.

\item RTDP\_EpsilonGreedy: A RTDP variant incorporating an epsilon-greedy exploration strategy.

\item Value\_Iteration\_UCB: An improved algorithm that introduces an exploration bonus under the value iteration framework\cite{tarbouriech2021stochastic}.

\end{itemize}

The convergence threshold for value updates was set to $\theta = 10^{-3}$, and the UCB exploration coefficient was fixed at $c=2$. The simulation consists of 100 runs, with each run comprising 300 episodes. In each episode, the AV starts from one of the designated origin states $s_{o}$ and proceeds until the goal state $s_{d}$. 
\begin{figure}[h!]
    \centering
\includegraphics[width=1\columnwidth]{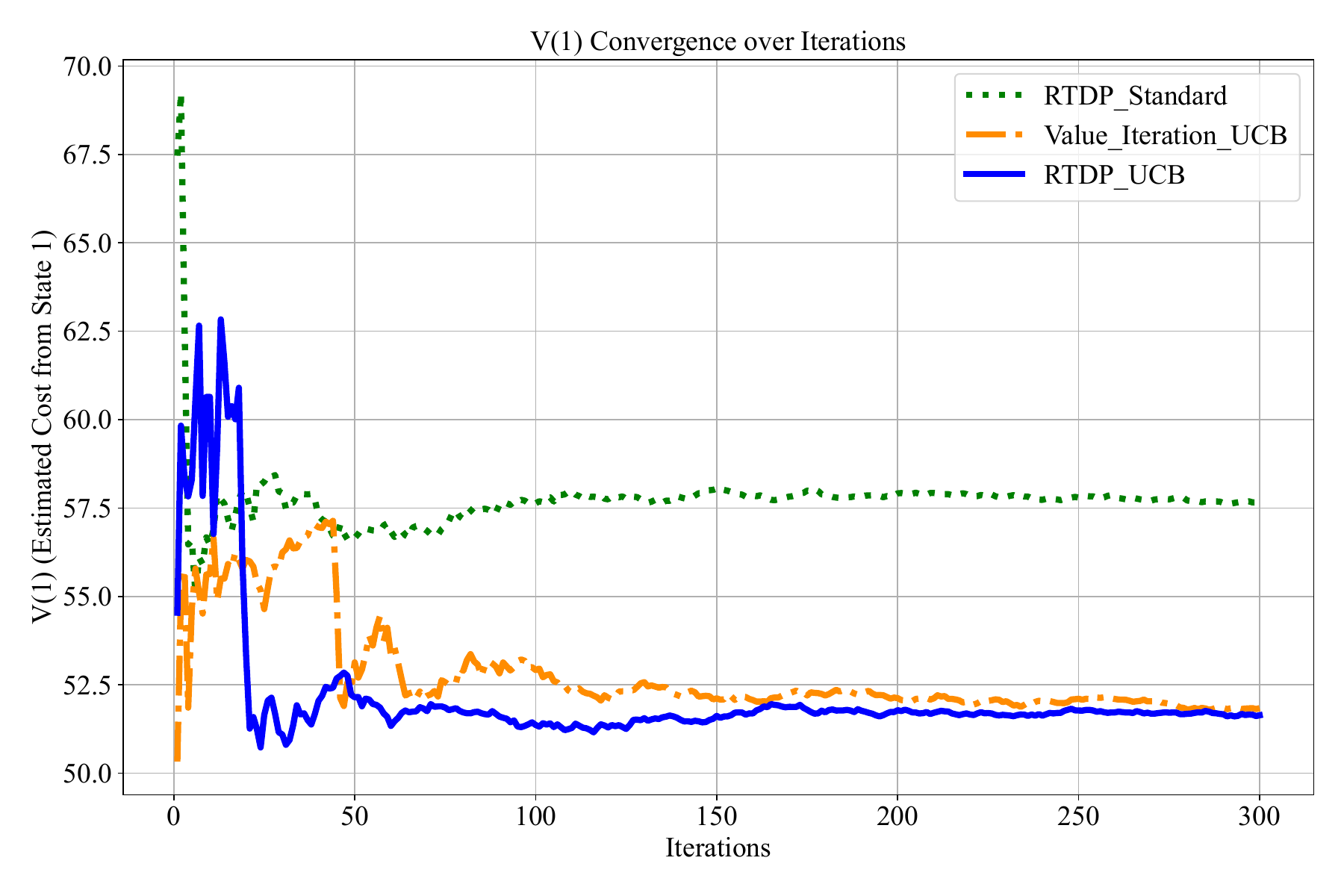}  
    \caption{Convergence of estimated state value $V(s_{o})$ over 300 iterations for four algorithms.}
    \label{fig:V1}
\end{figure}

Figure~\ref{fig:V1} illustrates the estimated state value at $s_{o}$ at the end of each iteration across all algorithms. As shown in Figure~\ref{fig:V1} , 
\begin{figure}[h!]
    \centering
    \includegraphics[width=1\columnwidth]{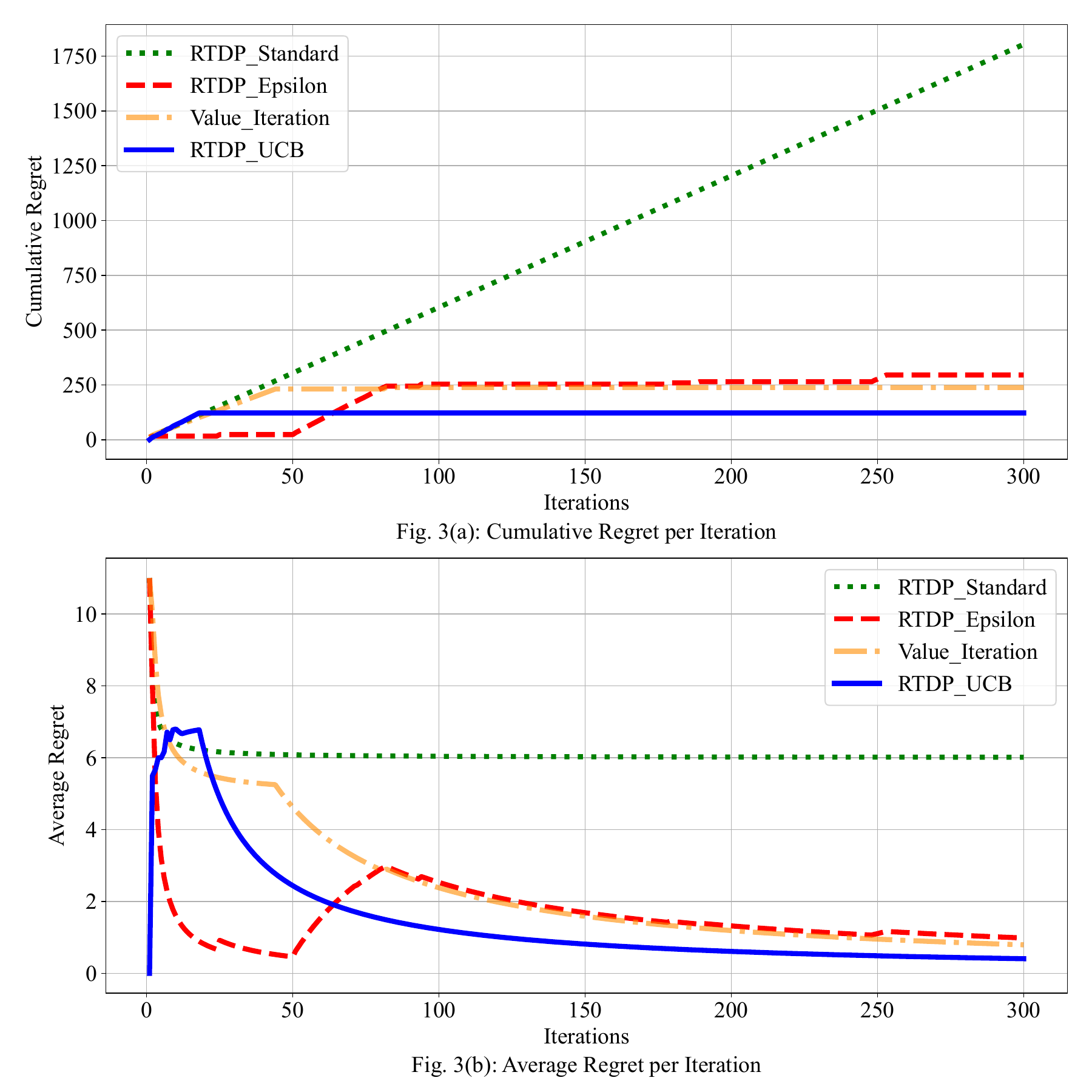}  
    \caption{Cumulative and average regret over 300 episodes for the four evaluated algorithms. }
    \label{fig:regret}
\end{figure}
the standard RTDP algorithm, due to its greedy decision-making strategy and lack of global value updates, is prone to getting trapped in locally optimal paths, particularly in environments with high stochasticity, thereby hindering effective convergence. In contrast, algorithms incorporating explicit exploration strategies—such as RTDP\_UCB and Value\_Iteration\_UCB—enable more comprehensive exploration of the state space during path planning and ensure that the final state values $V(1)$ converge to the true cumulative costs along the optimal paths.

The cumulative and average regret across all 300 episodes are presented in Figure~\ref{fig:regret} (a) and \ref{fig:regret}(b), respectively. The results demonstrate that incorporating exploration strategies, such as UCB or Epsilon-Greedy, within the RTDP framework, as well as integrating UCB into the standard Value Iteration framework, significantly enhances the convergence behavior of the algorithms. Among these algorithms, RTDP\_UCB exhibits the fastest convergence, highlighting the effectiveness of UCB-guided exploration in improving learning efficiency under stochastic environments. 

In terms of computational efficiency, Value\_Iteration\_UCB updates all state values in each iteration, offering stronger convergence guarantees but at a higher computational cost. In contrast, the asynchronous RTDP framework significantly reduces runtime. Average runtimes are reported in Table ~\ref{tab:runtime} . 

\begin{table}[htbp]
\centering
\caption{Comparison of algorithms: runtime, $V(1)$ estimate, and average regret}
\label{tab:runtime}
\begin{tabular*}{0.48\textwidth}{@{\extracolsep{\fill}}cccc}
\hline
\textbf{Algorithm} & \textbf{Time (s)} & \textbf{Est. $V(1)$} & \textbf{Avg. Regret} \\
\hline
RTDP\_Standard & 0.0340 & 57.63 & 6.01 \\
RTDP\_EpsilonGreedy & 0.0212 & 52.11 & 0.98 \\
Value\_Iteration\_UCB & 0.1480 & 51.84 & 0.79 \\
RTDP\_UCB & \textbf{0.0354 }& 51.64 & \textbf{0.41} \\
\hline
\end{tabular*}
\end{table}
\begin{figure}[h!]
    \centering
    \includegraphics[width=1\columnwidth]{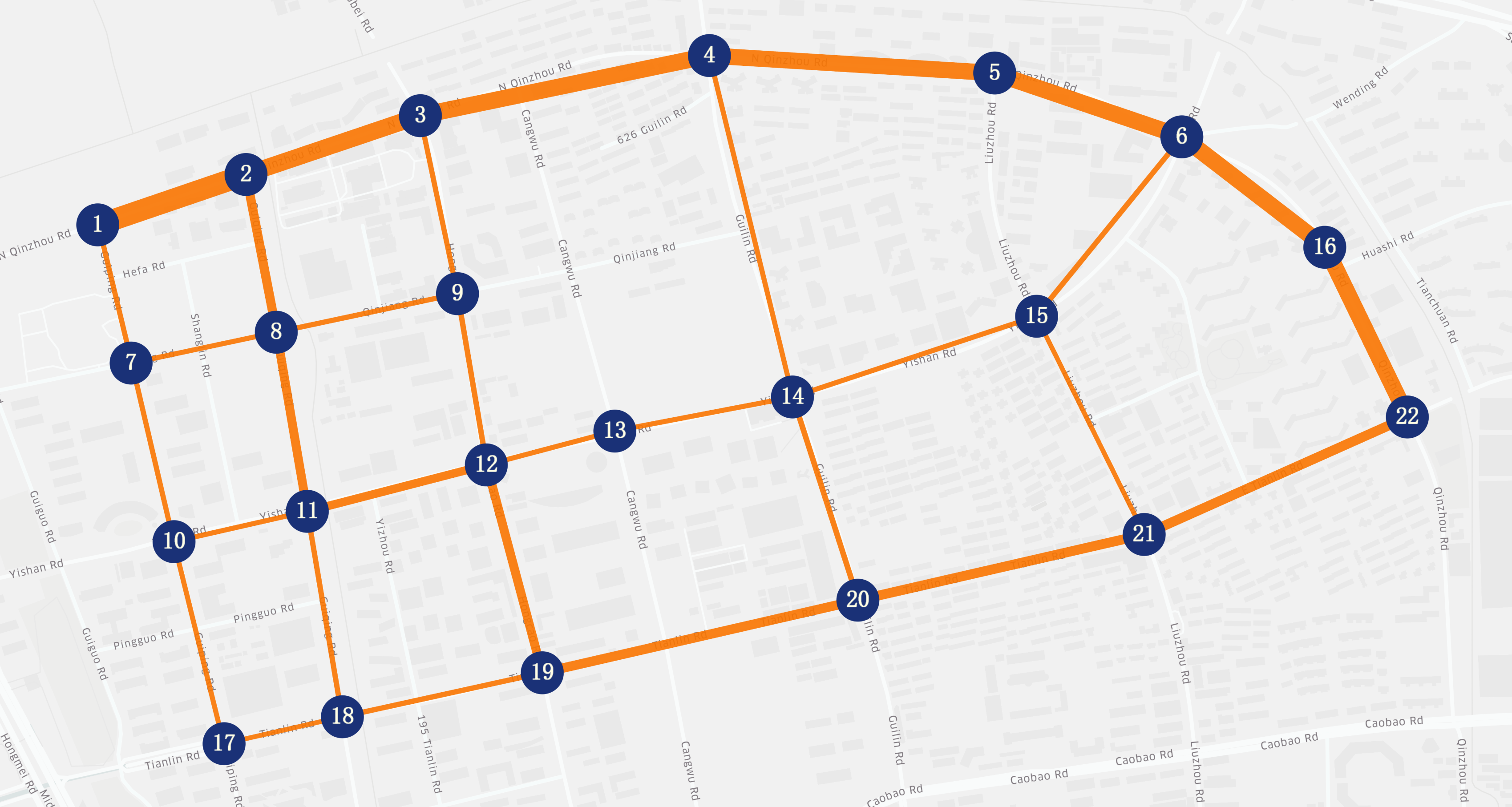}  
    \caption{Visualization of edge sampling during the simulation. The thickness of each edge represents the number of times it was sampled during the execution of the algorithm, indicating a certain degree of selection preference. }
    \label{fig:network}
\end{figure}

Although the algorithm's objective is to find the shortest path from $s_{o}$ to $s_{d}$, Figure \ref{fig:network} reveals that all states exhibit clear preferences for certain actions during the learning process. This indicates that the algorithm performs extensive policy learning across the entire state space and eventually converges to the optimal policy at every state. 

\section{CONCLUSIONS}
\label{05_conclusions}

In this paper, we propose an algorithm that integrates the asynchronous RTDP framework with a UCB exploration strategy to solve the routing problem of autonomous vehicles (AVs) in mixed-autonomy environments. By introducing optimism-based action selection into the RTDP process, the algorithm effectively balances exploration and exploitation, enabling efficient policy learning over large-scale state spaces and highly stochastic network environments. We derive and prove a path-level regret upper bound under the RTDP framework, confirming the sublinear growth of cumulative regret and the convergence guarantee of the algorithm. Furthermore, numerical experiments validated the practical efficiency of RTDP-UCB, showing faster convergence and lower average regret compared to benchmark algorithms. 

Future work may explore extending the RTDP-UCB framework to more general decision-making settings. First, it can be adapted to contextual bandit environments to incorporate observed contextual information into the exploration-exploitation trade-off. Second, the framework can be extended to partially observable Markov decision processes (POMDPs) to handle limited state observability. Third, extending the framework to consider origin-destination (OD) matrices would enable it to handle multiple concurrent routing tasks. 

\bibliographystyle{ieeetr}
\bibliography{ssp}

\end{document}